\documentclass[12pt]{iopart}
\usepackage{iopams}
\begin{document}

\title{On some sequences of functions and their applications in physics}

\author{Pablo G A Braz e Silva\dag
\footnote[3]{Partially supported by CAPES - BEX1901/99-0 -
Bras\' \i lia - Brasil}
and Andr\' es R R Papa\ddag
\footnote[4]{Also at International Centre for Theoretical Physics, Trieste,
Italy}}

\address{\dag\ Department of Mathematics and Statistics, University of New
Mexico, Albuquerque, NM 87131 USA}

\address{\ddag\ Instituto de Radioprote\c c\~ ao e Dosimetria IRD/CNEN\\
Av. Salvador Allende s/n, Recreio, 22780-160
Rio de Janeiro, RJ, Brazil}

\begin{abstract}
 We rigorously state and prove properties previously
used to ilustrate non-uniqueness of some generalized physical quantities.
This fact imposes peculiar constraints to the
mathematical functions that can be used to
represent such quantities.
\end{abstract}

\newtheorem{lemma}{Lemma}




\section{Introduction}

During the last decades many attempts
have been made in order to explain and formalize mathematically
a large set of phenomena that, at first glance, are not well described
by traditional theories. Among those attempts we can find the whole
theory of self-organized criticallity (that now seems to be over at least
in the criticallity part) and some statistical theories. Probably the most
popular in recent years has been the Dar\' oczy-like entropy
(Dar\' oczy 1970),
and its correponding
Statistical Mechanics formalism known as Tsallis entropy (Tsallis 1988):
\begin{equation}\label{eq1}
   S_q = k_B \frac{1 - \displaystyle \sum_{i=1}^{W} p_i^q}{q-1} ,
\end{equation}
where $W$ is the total number of configurations, $p_i$ are the associated
probabilities and $k_B$ is some suitable constant. It is based in a
special parameter $q$, which is adjusted to the phenomena the theory tries to
represent.
When $ q \rightarrow 1$ the Boltzmann entropy is obtained.
In Papa (1998), an inverse formalism was applied to find the limit when
$q \rightarrow 1$.
In this way, an infinite set of entropies(and consequentely
Statistics) of the same type was
found. The necessity of a deeper criteria to choose a formalism
to represent physical systems was pointed out.
In particular, it was shown that for the special
functional form of equation (\ref{eq1}),
integrating an infinite number of times the numerator and the denominator
separately leads
to the same result as differentiating once  ({\it i.e.}, to apply the replica trick or,
even simpler, L'H\^ opital's rule).
So, letting the number of integrations go to infinity, we achieve a limit that is independent of the parameter $q$, without considering the limit of $q$ approaching some specific value.
The aim of this paper is to state and prove a general Lemma providing
general conditions to apply this technique.

\section{The lemma}

For a real continuous function
$f:\mathbb{R} \rightarrow \mathbb{R}$ and a point
$x_0\in \mathbb{R}$, we denote by
$\displaystyle I_n^f (x)$ the sequence of functions
defined inductively as follows:
\begin{eqnarray*}
   I_0^f (x) = f(x) \\
   \displaystyle I_{n+1}^f (x) = \int_{x_0}^x I_n^f (t) dt .
\end{eqnarray*}
To keep the notation as simple as possible, we do not write
explicitely the dependence of the sequences of functions above
on the point $x_0$. The point used as lower limit
of integration should be obvious in what follows.
We use the usual notation
$f^{(j)} (x)$ to denote the $j$-th derivative of $f(x)$.
With this basic notation, we can state the lemma:
\begin{lemma}
 Let $f:\mathbb{R} \rightarrow \mathbb{R}$,
 $g:\mathbb{R} \rightarrow \mathbb{R}$ be real functions, analytic
in a neighborhood $U$ of $x_0 \in \mathbb{R}$. Suppose $f^{(j)} (x_0 ) =
g^{(j)} (x_0 ) = 0$, for $j = 0,1,...,K$ and
$g^{(K+1)} (x_0 ) \neq 0$. Let
$V = [x_0 - R , x_0 +R] \subset U$, $R>0$, be an interval where
$f^{(j)} (x) \neq 0$ and
$g^{(j)} (x) \neq 0$, $\forall x\in V$, $x\neq x_0$, for
 $j = 0,1,...,K+1$.
Then,
\[
 \displaystyle  \lim_{n\rightarrow \infty}\frac{I_n^f (x)}{I_n^g (x)} =
  \frac{f^{(K+1)}(x_0)}{g^{(K+1)}(x_0)}
\]
uniformily in $V\backslash \{ x_0\}$, where the sequences of functions $I_n$ are constructed
using $x_0$ as the lower limit of integration.
\end{lemma}
\paragraph{Proof:}
  Since $f$ and $g$ are analytic in V, we have for $x\in V$
\begin{eqnarray*}
  \displaystyle f(x) = \sum_{j=0}^{\infty} \frac{f^{(j)} (x_0)}{j!} (x-x_0)^j \\
  \displaystyle g(x) = \sum_{j=0}^{\infty} \frac{g^{(j)} (x_0)}{j!} (x-x_0)^j .
\end{eqnarray*}
$f^{(j)} (x_0 ) = g^{(j)} (x_0 ) = 0$, for $j = 0,1,...,K$, implies
\begin{eqnarray*}
  \displaystyle f(x) = \sum_{j=K+1}^{\infty} \frac{f^{(j)} (x_0)}{j!} (x-x_0)^j \\
  \displaystyle g(x) = \sum_{j=K+1}^{\infty} \frac{g^{(j)} (x_0)}{j!} (x-x_0)^j .
\end{eqnarray*}
 Since both series converge uniformly, we can integrate term by
term to get
  \begin{eqnarray*}
     \displaystyle I_n^f (x) =
\displaystyle \sum_{j=K+1}^{\infty}
\displaystyle \frac{f^{(j)} (x_0)}{(j+n)!} (x-x_0)^{(j+n)} \\
   \displaystyle I_n^g (x) =
\sum_{j=K+1}^{\infty} \frac{g^{(j)} (x_0)}{(j+n)!} (x-x_0)^{(j+n)} \\
\end{eqnarray*}
and then, for $x \neq x_0$,
\begin{eqnarray*}
\displaystyle \fl \frac{I^f_n (x)}{I^g_n (x)}  =
\frac{\displaystyle \sum_{j=K+1}^{\infty}
\frac{f^{(j)} (x_0)}{(j+n)!} (x-x_0)^{(j+n)}}
{\displaystyle \sum_{j=K+1}^{\infty} \frac{g^{(j)} (x_0)}{(j+n)!} (x-x_0)^{(j+n)}} \\
 \fl =  \frac{f^{(K+1)}(x_0) + \displaystyle \frac{1}{k+2+n}\left[
f^{(K+2)}(x_0) (x-x_0) + \sum_{j=2}^{\infty}
\frac{f^{(K+1+j)}(x_0)}{(K+1+j+n)...(K+3+n)} (x-x_0)^j\right]}{g^{(K+1)}(x_0) +
\displaystyle \frac{1}{k+2+n}\left[
g^{(K+2)}(x_0) (x-x_0) + \sum_{j=2}^{\infty}
\frac{g^{(K+1+j)}(x_0)}{(K+1+j+n)...(K+3+n)} (x-x_0)^j\right]}.
\end{eqnarray*}
Since
\begin{eqnarray*}
\displaystyle \fl \sum_{j=2}^{\infty}
 \frac{f^{(K+1+j)}(x_0)}{(K+1+j+n)...(K+3+n)} (x-x_0)^j & \hspace{.2cm}\mbox{and} \hspace{.2cm} &
\displaystyle \sum_{j=2}^{\infty}
\frac{g^{(K+1+j)}(x_0)}{(K+1+j+n)...(K+3+n)} (x-x_0)^j
\end{eqnarray*}
converge
uniformily in $V$, there exists a constant $C>0$ such that for all
$x \in V$,
\begin{eqnarray*}
\left| f^{(K+2)}(x_0) (x-x_0) + \sum_{j=2}^{\infty}
\frac{f^{(K+1+j)}(x_0)}{(K+1+j+n)...(K+3+n)} (x-x_0)^j\right| \leq C \\
\left| g^{(K+2)}(x_0) (x-x_0) + \sum_{j=2}^{\infty}
\frac{g^{(K+1+j)}(x_0)}{(K+1+j+n)...(K+3+n)} (x-x_0)^j\right|\leq C
\end{eqnarray*}
and that implies
\[
  \displaystyle  \lim_{n\rightarrow \infty}\frac{I_n^f (x)}{I_n^g (x)} =
  \frac{f^{(K+1)}(x_0)}{g^{(K+1)}(x_0)}
\hspace{.3cm}, \hspace{.3cm}\forall x\in V \backslash \{ x_0\}. 
\]

\section{Conclusion}

We showed that under quite general conditions,
the physically relevant limit can be reached from a special generalized
functional form, without considering the limit to the value of the
parameter that recover this point through traditional methods
(replica trick or
L'H\^ opital's rule). This curious fact imposes very strong constraints
to functions used to generalize known physical quantities.
Extra conditions
are necessary to choose one of them.

\ack
 The authors wish to thank The International Centre
for Theoretical Physics, Trieste, Italy, since their collaboration
was made possible
during a visit to the centre.

\References

\item[] Dar\' oczy Z 1970 {\it Inf. Contro.} {\bf 16} 36

\item[] Papa A R R 1993 {\it J. Phys. A: Math. Gen.} {\bf 31},
5271-5276

\item[] Tsallis C 1988 {\it J. Stat. Phys.} {\bf 52} 479

\endrefs

\end{document}